\documentclass[12pt]{article}
\usepackage{amsfonts}
\usepackage{amsmath}
\usepackage{amssymb}
\usepackage{amscd}
\usepackage{color}
\usepackage{graphics}
\usepackage{graphicx}
\usepackage{epsfig}
\begin{document}
\newtheorem{Theoreme}{Th\'eor\`eme}[section]
\newtheorem{Theorem}{Theorem}[section]
\newtheorem{Th}{Th\'eor\`eme}[section]
\newtheorem{De}[Th]{D\'efinition}
\newtheorem{Pro}[Th]{Proposition}
\newtheorem{Lemma}[Theorem]{Lemma}
\newtheorem{Proposition}[Theorem]{Proposition}
\newtheorem{Lemme}[Theoreme]{Lemme}
\newtheorem{Corollaire}[Theoreme]{Corollaire}
\newtheorem{Consequence}[Theoreme]{Cons\'equence}
\newtheorem{Remarque1}[Theoreme]{Remarque}
\newtheorem{Convention}[Theoreme]{{\sc Convention}}
\newtheorem{PP}[Theoreme]{Propri\'et\'es}
\newtheorem{Conclusion}[Theoreme]{Conclusion}
\newtheorem{Ex}[Theoreme]{Exemple}
\newtheorem{Definition}[Theoreme]{D\'efinition}
\newtheorem{Remark1}[Theorem]{Remark}
\newtheorem{Not}[Theoreme]{Notation}
\newtheorem{Nota}[Theorem]{Notation}
\newtheorem{Propo}[Theorem]{Proposition}
\newtheorem{exercice1}[Th]{Lemme-Confii au lecteur}
\newtheorem{Corollary}[Theorem]{Corollary}
\newtheorem{PPtes}[Th]{Propri\'et\'es}
\newtheorem{Def}[Theorem]{Definition}
\newtheorem{Defi}[Theorem]{Definition}
\newtheorem{Example1}[Theorem]{Example}
\newenvironment{Proof}{\medbreak{\noindent\bf Proof }}{~{\hfill $\bullet$\bigbreak}}

\newenvironment{Demonstration}{\medbreak{\noindent\bf D\'emonstration
 }}{~{\hskip 3pt$\bullet$\bigbreak}} 

\newenvironment{Remarque}{\begin{Remarque1}\em}{\end{Remarque1}} 
\newenvironment{Remark}{\begin{Remark1}\em}{\end{Remark1}}
\newenvironment{Exemple}{\begin{Ex}\em}{~{\hskip
3pt$\bullet$}\end{Ex}} 
\renewenvironment{abstract}{\medbreak{\noindent\bf Abstract :
 }}{\bigbreak}
\newenvironment{Notation}{\begin{Not}\em}{\end{Not}}
\newenvironment{Notation1}{\begin{Nota}\em}{\end{Nota}}
\newenvironment{Example}{\begin{Example1}\em}{~{\hskip
3pt$\bullet$}\end{Example1}} 

\newenvironment{Remarques}{\begin{Remarque1}\em \ \\* }{\end{Remarque1}}

\renewcommand{\Re}{{\cal R}}
\newcommand{\Dim}{{\rm Dim\,}}
\renewcommand{\Im}{{\cal F}}
\newcommand{\finpreuve}{~{\hskip 3pt$\bullet$\bigbreak}}
\newcommand{\hp}{\hskip 3pt}
\newcommand{\hph}{\hskip 8pt}
\newcommand{\hphh}{\hskip 15pt}
\newcommand{\vp}{\vskip 3pt}
\newcommand{\vpv}{\vskip 15pt}
\newcommand{\IP}{{\mathbb{IP}}}
\newcommand{\rd}{{\mathbb R}^2}
\newcommand{\R}{{\mathbb R}}

\newcommand{\Hyper}{{\mathbb H}}
\newcommand{\Int}{{\mathbb I}}
\newcommand{\Boule}{{\mathbb B}(0,1)}
\newcommand{\Cantor}{{\mathbb K}}
\newcommand{\dist}{{\mbox{\em dist}}}
\newcommand{\K}{{\mathbb K}}
\newcommand{\B}{{\mathbb B}}
\newcommand{\Ho}{{\mathbb H}}
\newcommand{\Nat}{{\mathbb N}}
\newcommand{\N}{{\mathbb N}}
\renewcommand{\P}{{\mathbb P}}
\newcommand{\Esp}{{\mathbb E}}
\newcommand{\Complex}{{\mathbb C}}
\newcommand{\Ha}{{\cal H}}
\newcommand{\Harm}{{\bold H}}
\newcommand{\Lcal}{{\cal L}}
\newcommand{\ds}{\displaystyle}
\newcommand{\un}{\bold 1}
\newcommand{\Cone}{C(x,r,\varepsilon ,\Phi)}
\newcommand{\Cn}{C(x,2^{-n},\varepsilon ,\Phi )}
\newcommand{\Tranche}{W(x,r,\varepsilon,\Phi)}
\newcommand{\Wn}{W(x,2^{-n},\varepsilon,\Phi)}
\newcommand{\WFn}{W(x,2^{-n},\varepsilon,\Phi)\cap F}
\newcommand{\ovec}{\overrightarrow}
\newcommand{\red}{{\bold R}}
\newcommand{\dimH}{\dim_{\Ha}}
\newcommand{\diam}{\mbox{diam}}
\newcommand{\diamit}{\mbox{\em diam}}
\newcommand{\para}{\vskip 2mm}
\newcommand{\cod}{\stackrel{\mbox{\tiny cod}}{\sim}}
\newcommand{\cardit}{\mbox{\em card}}
\newcommand{\card}{\mbox{card}}
\newcommand{\Sphere}{{\mathbb S}_d}
\newcommand{\distit}{\mbox{\em dist}}
\newcommand{\Tri}{{\cal P}}
\newcommand{\LL}{{\mathcal L}}
\newcommand{\infess}{\mbox{inf\,ess}}
\newcommand{\supess}{\mbox{sup\,ess}}
\newcommand{\I}{\vert}

\definecolor{darkblue}{rgb}{0,0,.5}
\def\u{\underline }
\def\o{\overline}
\def\h{\hskip 3pt}
\def\hh{\hskip 8pt}
\def\hhh{\hskip 15pt}
\def\v{\vskip 8pt}
\def\vv{\vskip 15pt}
\font\courrier=cmr12
\font\grand=cmbxti10
\font\large=cmbx12
\font\largeplus=cmr17
\font\small=cmbx8
\font\nor=cmbxti10
\font\smaller=cmr8
\font\smallo=cmbxti10
\openup 0.3mm
%
%
%

\title{\color{blue} On Brownian flights}
\author{Athanasios BATAKIS\footnote{MAPMO} , Pierre LEVITZ\footnote{LPMC, Ecole Polytechnique} and Michel ZINSMEISTER{$^*$}}
\date{}
\maketitle
\begin{abstract}Let K be a compact subset of $\R^n$. We choose at random with uniform law a point at distance$ \varepsilon$ of K and start a Brownian motion
(BM) from this point. We study the probability that this BM hits K for the first time at a distance $\geq r$ from the starting point.
\end{abstract}
 \footnotetext{This work is partially financed by the ANR MIPOMODIM}

\section{Introduction and motivation.}
Porous materials, concentrated colloidal suspensions or physiological organs such that lung or kidney are systems developping large specific surfaces with a rich variety of shapes that influence the diffusive dynamics of Brownian particles. A typical example is the diffusion of water molecules in a colloidal suspension. NMR relaxation allows to measure the statistics of the flights of these molecules over long colloidal shapes such as proteins or DNA chains. It is thus tempting to rely these statistics to the geometry of the molecules, the goal being to probe shapes using this method. An ideal (but far-reaching) objective would be to make up a DNA-test for example using NMR relaxation.\\
This program has been developped in \cite{GKLSZ} where various kind of simulations or experiments have shown remarkable commun properties.\\
All the simulations measure the statistics of the same random phenomenon: an irregular curve or surface is implemented, consisting of a union of a large (but finite) number of equal affine pieces. Such a piece is chosen at random with uniform distribution and a random walker is started at some small but fixed distance from this piece inside the complement of the surface. One is interested in the law of the variable $X=$ the length of the flight, i.e. the distance between the starting point and the first hitting point on the surface of the random walker. Whatever shape the surface shows, the experiment shows the same behaviour 
\begin{equation}\label{one}
\P\left(X>r\right)\sim_{\infty}r^{d_{e}-d-2}\end{equation}
where $d_{e}$ is the dimension of the ambiant space and $d$ is the Minkowski dimension of the surface.\\
A first mathematical explanation of this behavior was outlined in \cite{GKLSZ}. The main purpose of the present paper is to give first a rigorous statement of this result and to prove it with minimal assumptions so that all the cases of the simulations are covered. This will be the content of the second paragraph.\\
The third section concerns an alternative approach of the result in a special but important case: a quasiconformal perturbation of the line in 2D. The paper being dedicated to Fred Gehring is one of the reasons of this section, but not the only one: indeed this alternative proof is neat and instructive.\\
There is a case of particular importance for physics and particularly polymer physics. It is the case of a curve being a self-avoiding walk, since this has been shown to be a good model for polymers. In 2D there are two evidences that the result should remain in this case : first intensive simulations performed by P.Levitz and secondly some computations by Duplantier using conformal field theory. We present in the last section these simulations and sho how, despite the fact that this case is not covered by results of section 3, one can use results about SLE to prove the result in this situation.\\
We would like to conclude this introduction with two remarks:\\
1) The result may be surprising since it involves only the Minkowski dimension while the experiment suggests that harmonic measure is involved (and it is!) and consequently that the result should depend on some multifractal property of the harmonic measure. This is not the case because of the law of the choice of the starting point. A completely different behavior would occur if instead of uniform law we had chosen harmonic measure and this case is extremely interesting since essentially it models the second flight: this case will be considered in a forecoming paper.\\
2) The paper \cite{GKLSZ} has been written while the second author was hosted by the laboratory of physics of condensed matter at Ecole Polytechnique whose members he thanks for their warm hospitality. This paper is an example of a successfull pluridisciplinar research and we would like to emphasize the fact that firstly problems coming from physics (especially polymer physics) are extremely rich mathematically speaking and that, secondly, modern Function Theory as we herited from great mathematicians as Fred Gehring is an extremely efficient tool to attack these challenging problems.\\
Our goal will be achieved if the reader gets convinced of this last statement after having read this paper. 

\section{Geometric considerations}
For convenience all domains in the next two sections of this paper will be assumed to have compact boundary.
An  estimate like (\ref{one}) cannot be true for every domain. Some geometric conditions are needed : one of them is that Minkowski dimension of the boundary exists, i.e. that 
$d=\lim_{\epsilon\to 0}\frac{\log\#N_{\epsilon}}{-\log\epsilon}$ exists, where $N_{\epsilon}$ is the minimal number of cubes of size $\epsilon$ needed to cover the boundary.

As we shall see later Whitney decomposition is a central tool in our proof; a second condition we have to impose is that the number of cubes of the Whitney decomposition at distance $r$  must also be comparable to  $r^{-d}$. A sufficient condition for this is given by the first property of NTA domains (\cite{JK}) called the ``corkscrew'' condition:
\begin{Definition} We say that a domain $\Omega$ satisfies the ``corkscrew'' condition if there exists $r_0>0$ and a constant $c>0$ such that for all $r<r_0$ and any $x\in \partial\Omega$ there exists $y\in\Omega$ such that
$cr<\dist(x,y)<r$ and $\dist(y,\partial\Omega)>cr$.
\end{Definition}

\begin{Proposition}
Under the corkscrew condition the number of Whitney cubes with size that intersect the level surface $\Gamma_r=\{x\in\Omega\,; \,\dist(x,\partial\Omega)=r\}$ is comparable to 
the minimal of cubes of size $r$ needed to cover the boundary.
\end{Proposition}
The straightforward proof is left to the reader. 

\section{The case of open sets in $\R^n$, $n\ge 3$}
Although  our approach can be easily adapted to open sets in the plane, we present it for $\R^n, \,n\geq 3$ for two reasons: First, Green function and related formulas being different we try to avoid writing everything twice. Secondly, a different approach is proposed in section \ref{quasiconf} for open sets in the complex plane, using the quasi-conformal theory.
\begin{Notation}
Given an open set $\Omega$, a point $x\in\Omega$ and a set $A\subset \bar\Omega$  we denote $\P(x\hookrightarrow_{\Omega} A)$ the probability that Brownian motion started at x touches $A$ before leaving $\Omega$. The Green function of $\Omega$ will be denoted $G_{\Omega}$ and let $G_n$ be the Green function of $\R^n$. The ball of center $x$ and radius $r$ is denoted by $\B(x,r)$ and the distance from $x$ to $\partial\Omega$ is denoted by $d_x$. We say that two quantities $A,B$ are ``comparable" (we denote $A\sim B$) if the exists a constant $c$ such that $\frac1c A\le B\le cA$. The harmonic measure of a set $F$ (usually a subset of $\partial\Omega$) at a point $x\in\Omega$ is denoted by $\omega_{\Omega}(x,F)$ or $\omega(x,F,\Omega)$.
\end{Notation}

\subsection{Green estimates}
\begin{Proposition}\label{Green}
Let $\Omega$ be a domain in $\R^n$, $n>2$ and $x,y\in\Omega$. There exists a universal constant $C$ depending only on $n$ such that for all $\ell\le \frac12$ $$\P\Big(x\hookrightarrow_{\Omega}\B(y,\ell d_y)\Big)\le C\left(\frac{d_y}{d_x}\right)^{n-2}\P\Big(y\hookrightarrow_{\Omega}\B(x,\ell d_x)\Big)$$
\end{Proposition}
\begin{Proof}
First, notice that by the maximum principle,  $G_{\B(x,d_x)}(x,z)\le G_{\Omega}(x,z)\le G_n(x,z)$ for all $z\in\B(x,d_x/2)$. It is then easy to check that there exists a constant $c=c(n)$ such that $c^{-1}G_n(x,z)\le G_{\Omega}(x,z)\le G_n(x,z)$ for all $z\in\B(x,\ell d_x)$. 

Therefore the function that assigns $s\mapsto (\ell d_y)^{n-2}G_n(s,y)$ is harmonic in $\Omega\setminus \B(y,\ell d_y)$, tends to $0$ at $\partial\Omega$ and takes values between $c^{-1}$ and $1$ on $\partial\B(y,\ell d_y)$.
The probability $\P\Big(x \hookrightarrow_{\Omega} \B(y,\ell d_y)\Big)$ is then comparable to $(\ell d_y)^{n-2}G_n(x,y)$.  In a similar way $\P\Big(y \hookrightarrow_{\Omega} \B(x,\ell d_x)\Big)$ is equivalent to $(\ell d_x)^{n-2}G_n(y,x)$ and the proof is complete since $G_n$ is symmetric.
\end{Proof}

We suppose from now on that the domain $\Omega$ is uniformly ``fat", i.e. that there exists a constant $c>0$ such that  for any $x\in\Omega$ and $r\le 1$, we have 
\begin{equation}\label{capacitycondition}
\mbox{cap}_{\B(x,2r)}\Big(\B(x,r)\cap \partial\Omega\Big)\ge c\, \mbox{cap}_{\B(x,2r)}\Big(\B(x,r)\Big)
\end{equation} 
In particular, this implies that there exists a uniform lower bound $L>0$ of the probability that Brownian motion started at $x$ hits the  boundary of $\Omega$ before leaving $\B(x,2d_x)$ (see \cite{Ancona5}, lemma 5). This notion (also called ``uniform capacity density condition") has previously been introduced in various contexts, cf. \cite {JW}, \cite{Ancona5}, \cite{hkm}.

\subsection{Main Results}
We consider a Whitney decomposition of $\Omega$ in {\bf dyadic} cubes $Q$ satisfying 
$c_1|Q|\le d(Q,\partial\Omega)\le c_2|Q|,$
where $c_1<1<c_2$ are positive constants (powers of $2$, for convenience) depending on $n$. 
For $t>0$ we note ${\mathcal Q}_{t}$ the subcollection of cubes of the Whitney decomposition that intersect the level surface $\Gamma_{t}=\{x\in\Omega\;;\; d_x=t\}$. 
\begin{Theorem}\label{first} Take $\varepsilon<r$, fix a cube  $Q_r\in{\mathcal Q}_r$ and consider, for every $Q\in {\mathcal Q}_{\varepsilon}$, its center $x_Q\in Q$. Then  
$\displaystyle\sum_{\displaystyle Q\in{\mathcal Q}_{\varepsilon}}P\Big(x_Q\hookrightarrow_{\Omega}Q_r\Big)$ is equivalent to $\left(\frac{r}{\varepsilon}\right)^{n-2}$.
\end{Theorem}
\begin{Proof}
According to proposition \ref{Green} we have 
$$\frac1C\sum_{Q\in{\mathcal Q}_{\varepsilon}}P\Big(x_{Q_r}\hookrightarrow_{\Omega}Q\Big)\le \left(\frac{\varepsilon}{r}\right)^{n-2}\sum_{ Q\in{\mathcal Q}_{\varepsilon}}P\Big(x_Q\hookrightarrow_{\Omega}Q_r\Big) \le C \sum_{Q\in{\mathcal Q}_{\varepsilon}}P\Big(x_{Q_r}\hookrightarrow_{\Omega}Q\Big),$$
where $x_{Q_r}$ is the center of the cube $Q_r$. We now show that $\displaystyle \sum_{Q\in{\mathcal Q}_{\varepsilon}}P\Big(x_{Q_r}\hookrightarrow_{\Omega}Q\Big)$ is equivalent to the harmonic measure $\omega_{\Omega}(x_{Q_r},\partial\Omega)$ of $\partial\Omega$ at $x_{Q_r}$ (in $\Omega$) which equals $1$.
For this purpose we will use hypothesis (\ref{capacitycondition}) together with an easy  control of multiple coverings. 

Take any $Q\in{\mathcal Q}_{\varepsilon}$ and consider the cube $3c_2 Q$ of same center but of $3c_2$ times the sidelength of $Q$  ($c_2$ being the constant of the Whitney decomposition). By condition (\ref{capacitycondition}) the probability for Brownian motion started anywhere in  $Q$ to exit $\Omega$ before exiting $3c_2 Q$ is bounded below by a positive constant $c$. Hence, the harmonic measure $\omega_{\Omega}(x_{Q_r},3c_2 Q\cap\partial\Omega)$ of  $3c_2 Q\cap\partial\Omega$ at $x_{Q_r}$ (in $\Omega$)  is greater than $c P\Big(x_{Q_r}\hookrightarrow_{\Omega}Q\Big)$. Summing over all cubes  $Q\in{\mathcal Q}_{\varepsilon}$ we get 
$$\displaystyle c\sum_{Q\in{\mathcal Q}_{\varepsilon}}P\Big(x_{Q_r}\hookrightarrow_{\Omega}Q\Big)\le \displaystyle \sum_{Q\in{\mathcal Q}_{\varepsilon}}\omega_{\Omega}\Big(x_{Q_r},3c_2Q\Big)$$

On the other hand, every $x\in\partial\Omega$ can only belong to a finite number of cubes $3c_2Q$; this proves that $\displaystyle \sum_{Q\in{\mathcal Q}_{\varepsilon}}\omega_{\Omega}\Big(x_{Q_r},3c_2Q\Big)$ is comparable to the harmonic measure of the boundary $\partial\Omega$  of $\Omega$ and therefore $\displaystyle \sum_{Q\in{\mathcal Q}_{\varepsilon}}P\Big(x_{Q_r}\hookrightarrow_{\Omega}Q\Big)$ has an upper bound depending only on $n$ and on condition's (\ref{capacitycondition}) constant. The lower bound is trivial by the ``fatness'' condition.
\end{Proof}

\begin{Theorem}\label{main1}
Choose $Q$ at random with uniform law in ${\mathcal Q}_{\varepsilon}$. The probability for a Brownian motion started at any point $x$ of $Q$ to hit $\Gamma_{r}$ before exiting $\Omega$ is comparable to $\displaystyle \frac{\#{\mathcal Q}_{r}}{\#{\mathcal Q}_{\varepsilon}}\left(\frac{r}{\varepsilon}\right)^{n-2}.$
\end{Theorem}

\begin{Proof} First observe that by the Harnack principle we can choose $x=x_Q$ the center of the cube $Q$. Since the cube $Q$ is arbitrarily chosen there are $\#{\mathcal Q}_{\varepsilon}$ possible choices. 

We consider now the cubes of ${\mathcal Q}_{r}$ and we define $S_r$ as the part of the boundary of $\cup_{Q\in {\mathcal Q}_{r}}Q$ separating the set $\displaystyle \Gamma_r=\{x\in\Omega \mbox{ such that } d_x=r\}$ and $\partial \Omega$. We say that a cube in ${\mathcal Q}_{r}$ has a seashore if part of its boundary is also part of $S_r$, and , in this case, this access to the ``sea" is at least a square whose diameter is $\geq$ to a constant depending only on $n$ times the size of the cube. We then consider the open set $U$ consisting in the union of the components of $\displaystyle\Omega\setminus \bigcup_{Q\in {\mathcal Q}_{r}}\overline{Q}$. We denote by $\tilde{S_{r}}$ the boundary of this open set: we are interested in the probability that  Brownian motion started at any point $x$ of $Q$ hits $\tilde{S_{r}}$ before exiting $\Omega$. Denote by $V$ the component of $U$ containing $Q$. Let  $O$ be a cube in ${\mathcal Q}_{r}$ having a seashore. Each one of its sides  touching $ \tilde{S_{r}}$  contains a dyadic square $R$ of $c_1/8$ times the size of $O$, $c_1$ being the constant in Whitney decomposition, such that $R\subset \tilde{S_{r}}$.   Let $O_L$ be a cube of the same center as $O$ but $1+c_1/2$ its size.We consider the dyadique cube $R'$ contained in $O_L$ of size $c_1/8$ times the size of $O$ vertically above $R$ and at distance $c_1/8$ from $O$ .

We use the Boundary Harnack Principle to prove that the probability that Brownian motion started at $x_{Q}$ leave $V$ through $R$ is comparable to 
$P\Big(x_Q\hookrightarrow_{V}\, R'\Big)$. 

By ``adapted cylinder" to a graph of a Lipschitz function we understand a vertical revolution cylinder of finite height centered on the graph. Let us remind the Boundary Harnack Principle : Let $u$ and $v$ be positive harmonic functions on  a Lipschitz domain vanishing on the graph between the adapted cylinder  (to a graph-component of the boundary) ${\mathcal C}$ and  the ``sub"-adapted cylinder $\tilde{\mathcal C}$ of the same center and revolution axis but of $\ell$ times the size, $\ell<1$.  If ${\mathcal C}'$  is the ``middle" cylinder of the same center and revolution axis but of $\frac{1+\ell}{2}$ times the size of ${\mathcal C}$. Then for all   $x\in\partial  {\mathcal C}'\cap V$
$$\frac{v(x)}{u(x)}\sim \frac{v(P)}{u(P)},$$
where $P$ is the intersection point of the revolution axis of the cylinder ${\mathcal C}'$ and of its boundary. The multiplicative constants in the equivalence relation depend on the ratio (revolution radius):(height) of ${\mathcal C}$ , on $\ell$, on the local Lipschitz norm of the boundary and on the dimension of the space $n$, see \cite{Ancona2}.

 Remark that $O_L$ only touches the neighbouring cubes of $R$. Clearly, $O_L\cap V$  is a Lipschitz domain (its boundary is composed of a finite union of squares). We can find a finite number of adapted cylinders ${\mathcal C}_i$ that do not touch the cubes $R'$  such that the ``sub"-adapted cylinders $\tilde{\mathcal C}_i$  of half their size cover the boundary of  $\partial R\cap \partial V$. Furthermore, we can ask these cylinders to touch the boundary of the cube $O_l$ of the same center as $O$ but of $1+c_1/100$ the size (see figure \ref{fig1}). Note that, since there is a finite number of configurations of the neighborhood of $O$ in the Whitney decomposition,  the number of adapted cylinders needed to this covering is bounded by a uniform constant $\kappa$ depending only on $n$.
We consider the harmonic functions $u(x)=\omega(x, R',V\setminus \bigcup_{\mbox{``wet'' sides of }O} R')$ and $v_i(x)=\omega(x,\tilde{\mathcal C}_i \cap \partial O_L, V)$  in the domain $V'=V\setminus \bigcup_{\mbox{``wet'' sides of } O} R'$. Clearly,
\begin{equation}\label{compare}
\omega(x, \partial O\cap\partial V,V')\le \sum_i v_i(x)\le \kappa \omega(x, \partial O\cap\partial V, V').
\end{equation}
\begin{figure}
\centerline{\epsfig{file=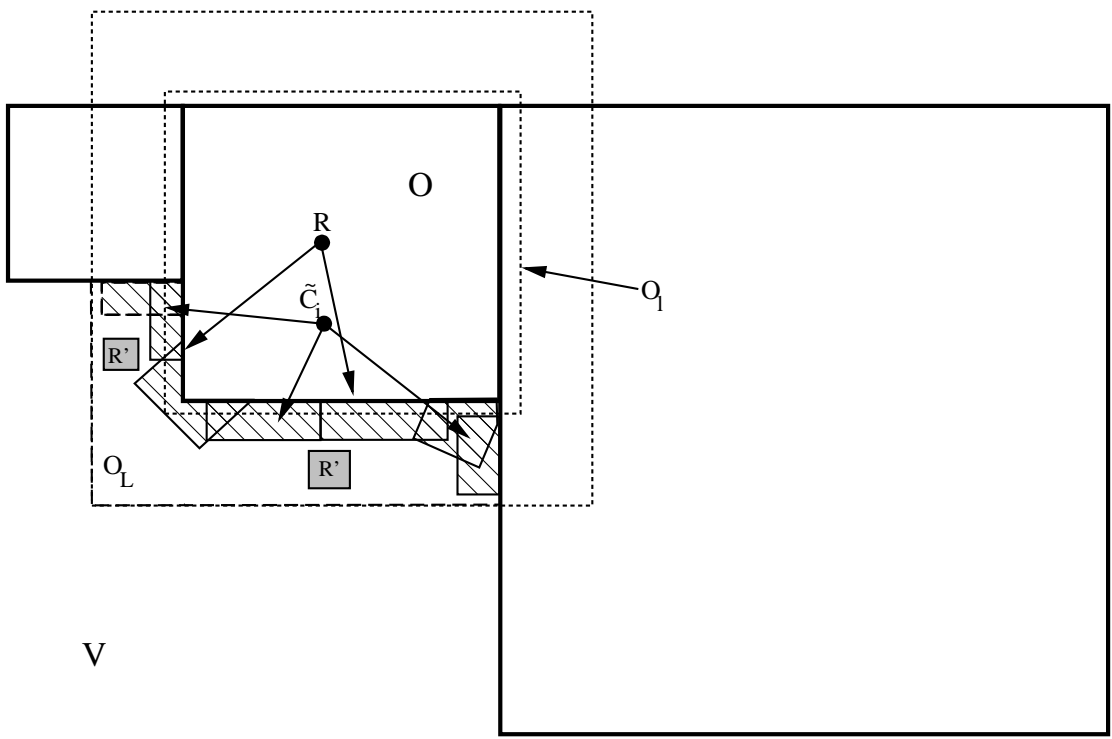, height=10cm}}
\caption{\nor A configuration in $\R^2$.}\label{fig1}
\end{figure}
We can apply the Boundary Harnack Principle to every one of these cylinders ${\mathcal C}_i$ for the functions $u$ and $v_i$. We get that for all $i$ and all $x\in\partial {\mathcal C}_i\cap V$
\begin{equation}\label{BHP}
\frac{v_i(x)}{u(x)}\sim \frac{v_i(P_i)}{u(P_i)},
\end{equation}
where $P_i$ is the intersection point of the revolution axis of the cylinder and of it's boundary. The boundary of our domain is the inner boundary of the cubes of the Whitney decomposition; hence we may restrain ourselves to a finite number of  configurations of the ``adapted cylinders" (up to a contractions-dilatations that do not affect constants)  and therefore the multiplicative constants in the equivalence relation will be finite in number and so uniformly bounded away from $0$ and infinity. Since the distance of  $\mbox{dist}(P_i,R')$  is equivalent to the distance $\mbox{dist}(P_i,\partial R)$ we can prove by standard arguments of harmonic analysis that $u(P_i)$ is bounded below by a constant depending only on dimension $n$. Using the Harnack principle and the fact that all adapted cylinders intersect $\partial R_l$ we get the existence of  a constant $c$ depending only on $n$ such that $v_i(P_k)\le c v_i(P_j)$ for all $i,j,k$.

After summing over $i$, taking in account equations (\ref{compare}) and (\ref{BHP}) and using standard Harnack inequalities we get
$$\omega(x, \partial R\cap\partial V', V')\le \omega(x, R',V') \le \kappa' \omega(x, \partial R\cap\partial V', V')$$ for all 
$x\in \partial R_l\cup\cup_i\partial {\mathcal C}_i\cap V$. By the maximum principle we thus obtain
$$\omega(x, \partial R\cap\partial V', V')\le \omega(x, R',V') \le \kappa' \omega(x, \partial R\cap\partial V', V')$$
for all $x\in Q\in {\mathcal Q}_{\varepsilon}$.

Theorem \ref{first} applies to the domain $V'$ where ${\mathcal Q}_r$ is replaced by the collection $R'$ coming from all ``wet''sides of all cubes with a seashore. We obtain that 
$$\displaystyle\sum_{\displaystyle Q\in{\mathcal Q}_{\varepsilon}, Q\subset V'}P\Big(x_Q\hookrightarrow_{ V'}R'\Big)\sim \left(\frac{r}{\varepsilon}\right)^{n-2}\omega_{V'}\left(x_{R'},\partial\Omega\cap V'\right)$$
a quantity bounded from below by some $\alpha>0$ by the fatness condition. The proof is completed by suming over all $R'$ and $V'$s and by noticing that up to a multiplicative factor depending only on the dimension, the number of cubes with a seashore is comparable to the number of cubes in $\mathcal {Q}_{r}$.
\end{Proof} 

\begin{Theorem}\label{main2}
Choose $Q$ at random with uniform law in ${\mathcal Q}_{\varepsilon}$. The probability for a Brownian motion started at any point $x$ of $Q$ to exit $\Omega$ at distance greater than $r$ from the starting point is comparable to $\displaystyle \frac{\#{\mathcal Q}_{r}}{\#{\mathcal Q}_{\varepsilon}}\left(\frac{r}{\varepsilon}\right)^{n-2}.$
\end{Theorem}
 \begin{Proof}
It depends on a comparison between the probability of ``cruising along the coast" $\partial\Omega$ and the probability to move at distance $r$ before coming back to the coast. The second is comparable to $\displaystyle \frac{\#{\mathcal Q}_{r}}{\#{\mathcal Q}_{\varepsilon}}\left(\frac{r}{\varepsilon}\right)^{n-2}$  according theorem \ref{main1} while the first is exponentially small. 
 
 Take $s>0$ and $x_Q\in Q\in{\mathcal Q}_{\varepsilon}$.  Consider the annuli centered at $x$ of inner radii $\ell s$ and outer radii $(\ell+1)s$ where $\ell=0,...,\left[\frac rs\right]$. Brownian motion started at $x$ and moving at distance $r$ from $x$ before exiting $\Omega$ must go through all these annuli. 
The probability of going through such an annulus while staying at distance at most $\frac s4$ from the boundary is bounded by a $p_0\in (0,1)$ by the ``fatness" hypothesis. To see this take any point $y$ in the middle of the annulus (i.e. at distance $\frac{\ell+1}2s$ from $x$) and consider the ball of center $y$ and radius $\frac s2$. If $d_y<\frac s4$, the probability to exit the ball without  touching $\partial\Omega$ is uniformly bounded away from $1$ by the ``fatness" hypothesis. This probability being greater than the probability of going through the annulus  we have the statement.
By the independence of the ``crossing annulli" events we get that the probability that Brownian motion  goes through all the annulli is smaller that $p_0^{\left[\frac rs\right] }$.

Let us now prove the following statement: ``there exist $0<c_1,c_2<1$ positive constants depending only on dimension and on the constant $L$ of the fatness condition such that for any $x\in\Omega$ there exist disjoint sets $K_1,K_2\subset\partial\Omega\cap\B(x,2d_x)$ verifying $\mbox{dist}(K_1,K_2)>c_1d_x$ and $\omega_{\Omega}(x,K_i)>c_2$ for $i=1,2$. Once more this is a consequence of the ``fatness" property. Cut the sphere $\partial\B(x,2d_x)$ in small equal normal polygons (spherical triangles in $\R^3$) and consider the intersections $L_i$ of the cones of summit $x$ and basis these polygons with the set $ \B(x,2d_x)\setminus\B(x,d_x)$ (see figure \ref{fig2}). Clearly, we can choose the polygons small enough (independently of $x$ and $d_x$) to have that
$\omega(x, L_i, \B(x,2d_x)\setminus L_i)<L/2^n$, where $L$ is the constant in the fatness condition. It is now clear that $\omega(x, L_i\cap\partial\Omega, \B(x,2d_x)\cap \Omega)<L/100$ and therefore there are two non-neighboring $L_i$'s having harmonic measure greater than $L/2^n\#L_i$, which proves the statement. 

\begin{figure}
\centerline{\includegraphics{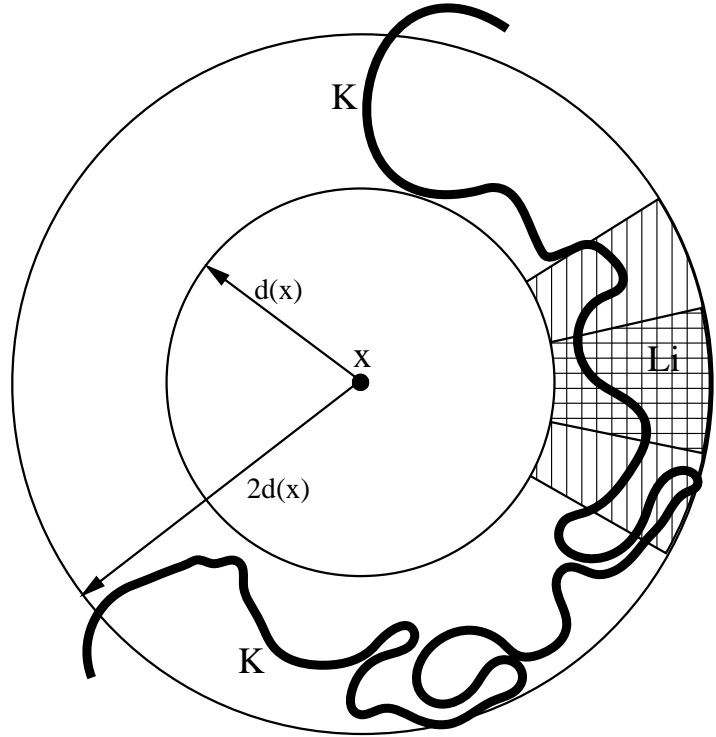}}
\caption{\nor Note that $d(x)$ is the distance of $x$ from the boundary $K=\partial\Omega$.}\label{fig2}
\end{figure}

The previous statement implies that once we reached  distance $r$ from the boundary of the domain the Brownian motion will revisit the boundary at distance comparable to $r$ from the starting point with probability greater than $c_2$. Putting together all the above we get the theorem. 
\end{Proof}
As a corollary we get the following.
\begin{Theorem}
If $\Omega$ satisfies
\begin{enumerate}\label{main}
\item The corkscrew condition
\item The fatness condition
\item and if $\partial\Omega$ has a Minkowski content, 
\end{enumerate}
then for every $\eta>0$  there exists a constant $c_{\eta,n}>0$ such that for all $r>\varepsilon>0$ if we choose $Q$ at random with uniform law in ${\mathcal Q}_{\varepsilon}$, the probability that a Brownian motion started at any point $x$ of $Q$ hits for the first time $\partial \Omega$ at distance greater than $r$ from the starting point $\P\left(X>r\right)$ verifies $$\displaystyle \frac{1}{c_{\eta,n}}\left(\frac{r}{\varepsilon}\right)^{n-d-2+\eta}\le\P\left(X>r\right)\le c_{\eta,n}\left(\frac{r}{\varepsilon}\right)^{n-d-2-\eta}.$$
\end{Theorem} 
\section{An alternative approach in the quasicircle perturbative case.}\label{quasiconf}
We present here a simple 2D case for which the proof of the main result is particularly simple using conformal mapping which preserves Brownian trajectories. The curve we consider will be a quasiconformal perturbation of the line. More precisely we consider a domain $\Omega=\varphi(\R^{2}_{+})$ where $\varphi:\R^{2}_{+}\to\mathbb C$ is holomorphic and such that
\begin{equation}
\sup_{x+iy\in \R^{2}_{+}}y\I\frac{\varphi''(x+iy)}{\varphi'(x+iy)}\I<1/2.\label{pert}
\end{equation}
It is known that under this hypothesis $\varphi$ is injective and has a quasiconformal extention to the whole plane. In particular $\Gamma=\varphi(\R)$ is a quasicircle close to a line. We will moreover assume that $\Gamma$ has a Minkowski dimension which we denote by $d$. By Koebe theorem, the quantity $y\I\varphi'(x+iy)\I$ is uniformly comparable to the distance from $\varphi(x+iy)$ to $\Gamma$. For $\alpha>0$ we then define $L_{\alpha}=\{x+iy\in \R^{2}_{+}\,;\,y\I\varphi'(x+iy)\I=\alpha\}$ and $\varphi(L_{\alpha})$ will serve as a substitute for the level set $\{\zeta\in\Omega\,;\,d(\zeta,\Gamma)=\alpha\}$.
\begin{Lemma}
If $\varphi$ satisfies ($\ref{pert}$) then $L_{\alpha}$ is the graph of a Lipschitz function $f_{\alpha}:\R\to\R$.
\end{Lemma}
\begin{Proof} We apply the implicit function theorem to the function\[F(x,y)=y\varphi'(x+iy)\overline{\varphi'}(x+iy)-\alpha^{2}.\]The computation gives\[\frac{\partial F}{\partial x}=2y\I\varphi'(x+iy)\I^{2}\Re(\psi(x+iy)),\frac{\partial F}{\partial y}=\I\varphi'(x+iy)\I^{2}(1+2y\Im(\psi(x+iy))\] where \[\psi(z)=\frac{\varphi''(z)}{\varphi'(z)}.\] The result follows because $\frac{\partial F}{\partial y}>0$ and by $(\ref{pert})$.\\
We now consider a portion of the curve with diameter $1$ and $0<\varepsilon<r$. We divide the portion of $\Omega$ between $\Gamma$ and $L_{\varepsilon}$ into pieces of diameter $\sim\varepsilon$. The preimages of these pieces are rough squares of sidelength $f_{\varepsilon}(x_{j})$ where $x_{j}$ is the left-hand point of the intersection with $\R$. By the results of the preceeding paragraph it suffices to compute the probability that a Brownian motion started at $\varphi(x_j+if_{\varepsilon}(x_{j}))$ will hit $\varphi(L_{r})$ before returning to $\Gamma$. By conformal invariance, this probability is comparable to $f_{\varepsilon}(x_{j})/f_{r}(x_{j})$ a quantity which is equivalent by Koebe to\[\frac{r}{\varepsilon}\frac{ \I \varphi'(x_{j}+if_{r}(x_{j})\I}{ \I \varphi'(x_{j}+if_{\varepsilon}(x_{j})\I}.\] On the other hand, by quasiconformality,\[(x_{j+1}-x_{j})\I \varphi'(x_{j}+if_{\varepsilon}(x_{j})\I\sim \varepsilon.\] Combining all the estimates we see that the probability we look for is comparable to\[ \varepsilon^{d}\left(\frac{r}{\varepsilon}\right)\sum_{j}\frac{ \I \varphi'(x_{j}+if_{r}(x_{j}))\I}{ \I \varphi'(x_{j}+if_{\varepsilon}(x_{j})\I}\sim\varepsilon^{d}\left(\frac{r}{\varepsilon}\right)\sum_{j}\frac{(x_{j+1}-x_{j})}{\varepsilon}\I \varphi'(x_{j}+if_{r}(x_{j}))\I\]and we see a Riemann sum appearing: we finally get as an estimate for the probability we seek\[\frac{\varepsilon^{d}}{r}\mathrm{length}(\varphi(L_{r}))\sim\frac{\varepsilon^{d}}{r}\frac{r}{r^{d}}\sim\left(\frac{r}{\varepsilon}\right)^{d},\]which is precisely what we wanted.
\end{Proof}
The preceeding proof has been presented because it is particurlarly simple, but of course the result is not optimal. The results of paragraph 2 remain true in dimension 2 and the proof requires only minor changes. In particular the result is true for all quasicircles; the difference with the case presented is that the topology of the level sets of the function distance to the boundary is more complicated in general.

\section{The self-avoiding walks case}

\begin{figure}
\centerline{\includegraphics[width=10cm,angle=-90]{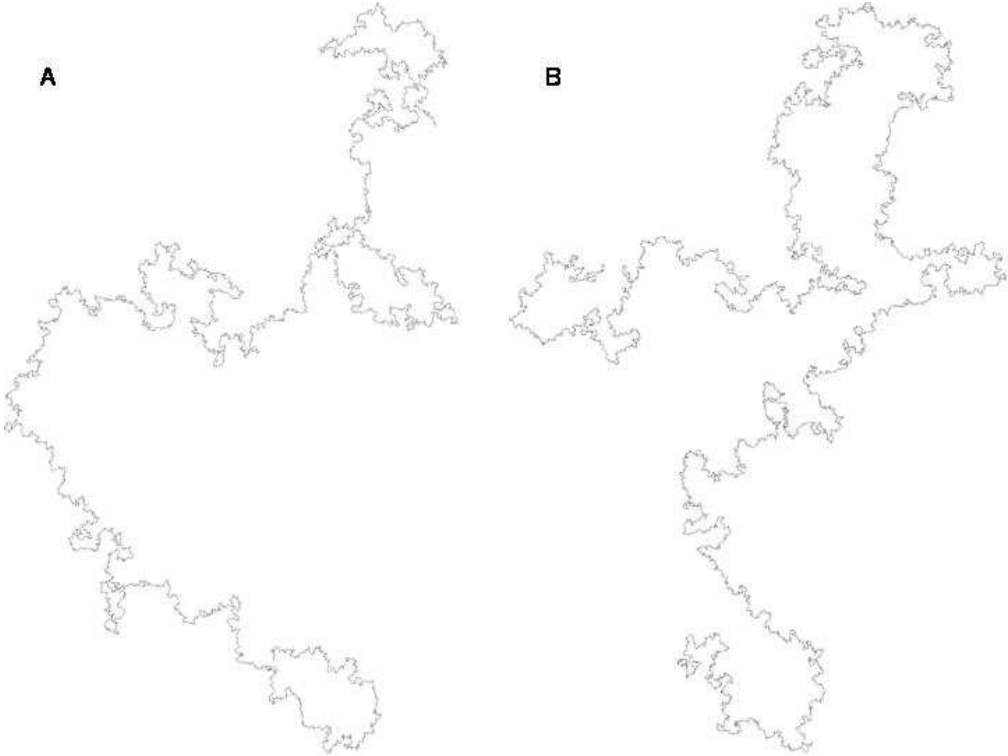}}
\caption{\nor Two examples of self avoiding walk in $2D$ generated by the pivot algorithm ~\cite{Kennedy}. The number of steps in each SAW is fixed at $10^{5}$.}
\label{plfig1}
\end{figure}

Self-avoiding walks (S.A.W.) (see ~\cite{madras} and ~\cite{PGG} for a definition) serve as a good model for polymers in physics. On the other hand it is strongly believed that in 2D self-avoiding walks is the same as $SLE_{8/3}$ (see \cite{RS} for a definition of $SLE_{\kappa}$). This conjecture, highly plausible, is comforted by the adequation between the computed dimension, which is $\displaystyle\frac43$ and the proved dimension for $SLE_{\kappa}$ curves , $1+\frac\kappa8$ (Beffara, ~\cite{Beffara}).
The following simulations can be seen as a new way of probing the adequation between SAW's and $SLE_{8/3}$:\\
In order to check that the statistics of flights over a self avoiding walk follow the expected law with $d=\frac43$, we have performed extended computer simulations. We have generated a set of self-avoiding walks on a square lattice using an implementation of the pivot algorithm described by Kennedy ~\cite{Kennedy} (see also \cite{Kennedy1}). The number of steps of the self avoiding walk is fixed at $10^{5}$. Two S.A.W. are shown in Fig.~\ref{plfig1}.

\begin{figure}[ht]
    \begin{center}  
      \includegraphics[scale=0.6,angle=-90]{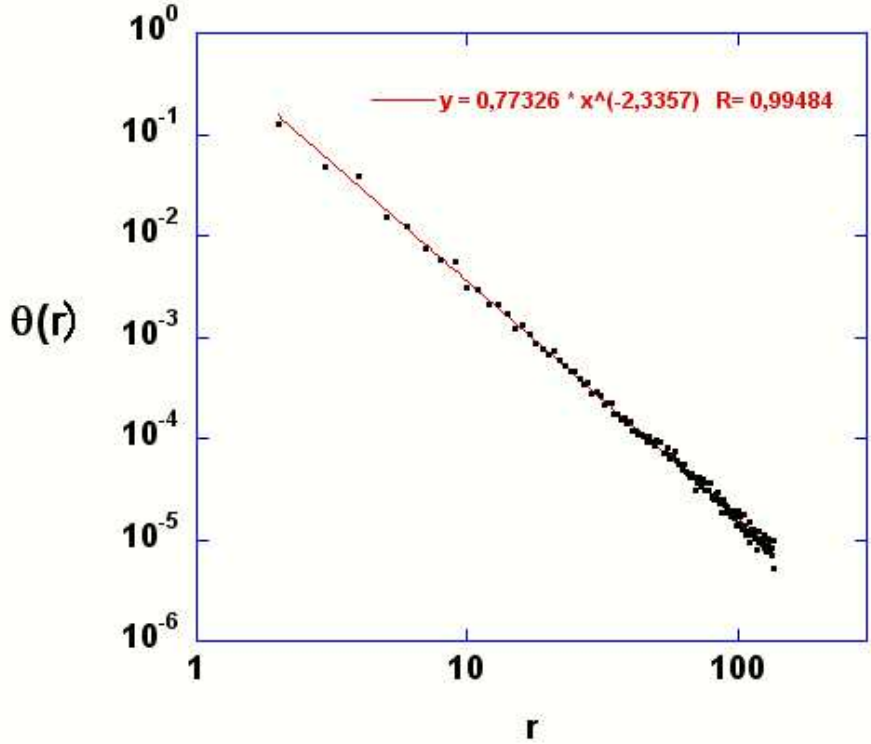}
    \caption{\nor Evolution of the probability density $\theta(r)$ that a particle
starting from close vicinity of SAW, returns to
the SAW, for the first time, after an end to end displacement found between $r$ and ${r+dr}$. The numerical estimation of the exponent $\beta$ is very close to $7/3$.}
\label{plfig2}
    \end{center}
\end{figure}

\begin{figure}[ht]
    \begin{center}  
      \includegraphics[scale=0.6,angle=-90]{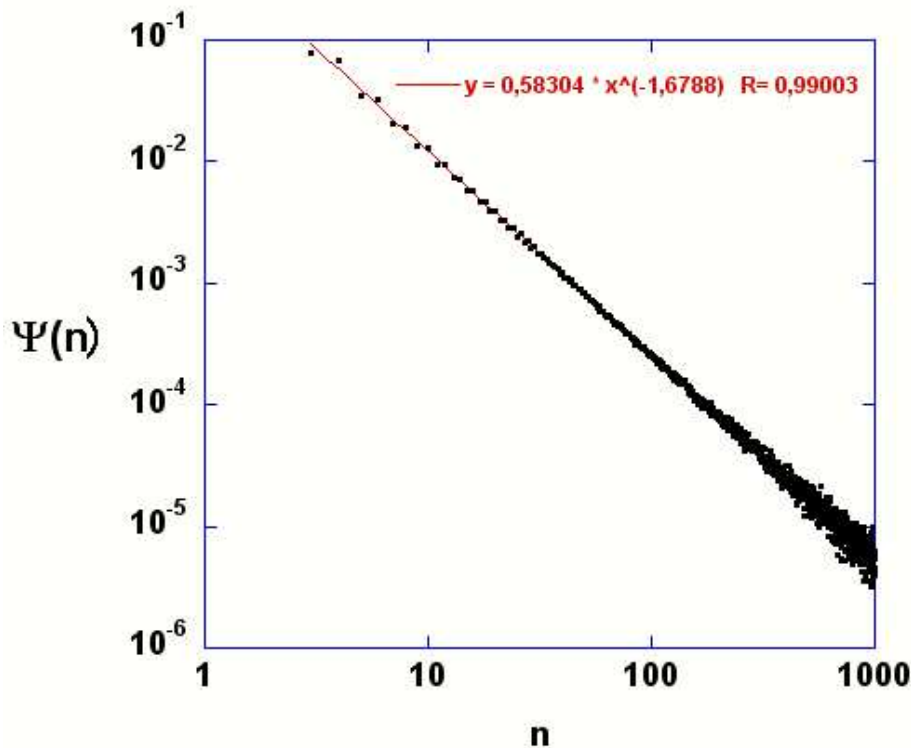}
    \caption{\nor Evolution of the probability density $\psi(n)$ that the particle
starting from a close vicinity of SAW and returning for the first time to
this SAW, has a total length displacement between $n$ and
$n+1$. The numerical estimation of the exponant $\alpha$ is very close to $10/6$.}
\label{plfig3}
    \end{center}
\end{figure}

We have checked by a box counting method that the mass fractal dimension of these curves are numerically founded around $1.33\pm0.005$. These values are very close to the expected value $4/3$. 
We have performed an on-lattice simulation analyzing the first passage statistics of a random walk starting in the close vicinity of the SAW and going back for the first time nearby the SAW. The numerical computations were performed on several configurations of SAW, using a statistical analysis over more than $2 10^9$ flights. Two probability density functions were computed. First, the probability density $\psi(n)$ that a flight has a total length equals to n. Second, the probability distribution of displacements $\theta(r)$. In order to limit edge effects, we have selected flights starting and ending on the same side of the S.A.W.. As shown Fig.~\ref{plfig2} and Fig.~\ref{plfig3}, we found that $\psi(n) \propto n^{-\alpha}$  and $\theta(r)\propto r^{-\beta}$.




It was shown in ~\cite{GKLSZ} that for a boundary of fractal dimension $d$ embedded in an  Euclidian space of Euclidean
dimension $d_e$, we should have

\begin{equation}
\label{eq-alpha}
\alpha = \frac{d - d_e + 4}{2} .
\end{equation}  
We are now in position to prove this estimate rigorously; this result is contained a paper to come.

Moreover as  $\theta(r)=-dP(r)/dr$, we get
\begin{equation}
\label{eq-beta}
\beta = {d - d_e + 3} .
\end{equation}

For $d=4/3$, we expect to have $\alpha=10/6$ and $\beta=7/3$. As shown in Figs.~\ref{plfig2}and ~\ref{plfig3}, numerical results provide a very good approximation of these above predictions.

After these convincing simulations, let us prove rigourously these asymptotics for $SLE_{\kappa}$ curves. First, combining results of Rohde-Schramm \cite{RS} and Beffara \cite{Beffara}, we see that the limsup in the definition of Minkowski upper-dimension for $SLE_{\kappa}$ is actually a limit, allowing asymptotic values for all values of $r$. Secondly, we observed in the last section that the result follows from the understanding of the number of Whitney cubes of given order. By a nice result of Bishop \cite{Bishop3} , it follows immediately that we get the right estimate if we allow the starting point to be chosen on both sides of the curve, which is actually the case in the above simulations.\\
If we always start from the same side, then not only Bishop's result does not apply but neither does it follow from the previous discussion because the corresponding domains do not have the corkscrew property. But, as Rohde and Schramm have proved, these domains are H\" older, meanning that the Riemman mapping from the upper half-plane onto them is H\"older continuous. This condition implies a weaker form of the the corkscrew condition which is sufficient to ensure the possibility to compute the Minkowski dimension of the $SLE_{\kappa}$ curve via the counting of Whitney cubes. As we have seen this is enough to prove the main result about statistics of flights.

The case of self affine curves is a little more delicate and will be treated in a forthcoming paper.

 It is to be noticed that, using quantum gravity arguments, Duplantier  also obtained the right exponents for all  $SLE_{\kappa}$(cf. \cite{Duplantier})

\bibliographystyle{alpha}
\bibliography{biblio}
\end{document}